\documentclass{amsart}
\usepackage{amsmath,amsfonts,amsthm,amssymb,color,graphicx,subfigure}
 \numberwithin{equation}{section}
 \theoremstyle{plain}
   \newtheorem{theorem}{Theorem}[section]
   \newtheorem{proposition}[theorem]{Proposition}
   \newtheorem{lemma}[theorem]{Lemma}

\theoremstyle{definition}
   \newtheorem{definition}[theorem]{Definition}
\usepackage{fullpage}
\author{Jessica Striker}
\title{The poset perspective on alternating sign matrices}
\address{School of Mathematics, University of Minnesota, Minneapolis, MN 55455}
\email{jessica@math.umn.edu}

\begin{document}
\maketitle

\begin{abstract}
Alternating sign matrices (ASMs) are square matrices with entries 0, 1, or $-1$ whose rows and columns sum to 1 and whose nonzero entries alternate in sign. We put ASMs into a larger context by studying the order ideals of subposets of a certain poset, proving that they are in bijection with a variety of interesting combinatorial objects, including ASMs, totally symmetric self--complementary plane partitions (TSSCPPs), Catalan objects, tournaments, semistandard Young tableaux, and totally symmetric plane partitions. 
We use this perspective to prove an expansion of the tournament generating function as a sum over TSSCPPs which is analogous to a known formula involving ASMs.
\end{abstract}

\section{Introduction}
Alternating sign matrices (ASMs) are simply defined as square matrices with entries 0, 1, or $-1$ whose rows and columns sum to 1 and alternate in sign, but have proved quite difficult to understand (and even count).
Totally symmetric self--complementary plane partitions (TSSCPPs) are 
plane partitions, each equal to its complement and invariant under all permutations of the coordinate axes.
TSSCPPs inside a $2n \times 2n \times 2n$ box are equinumerous with $n\times n$ ASMs, but no explicit bijection between these two sets of objects is known. In this paper we present a new perspective which sheds light on ASMs and TSSCPPs and brings us closer to constructing 
a explicit ASM--TSSCPP bijection.

\section{The tetrahedral poset}
\label{sec:intro}
Given an $n \times n$ ASM $A$, consider the following bijection to objects called monotone triangles of order~$n$~\cite{BRESSOUDBOOK}.
For each row of $A$ note which columns have a partial sum (from the top) of 1 in that row. Record the numbers of the columns in which this occurs 
in increasing order. This gives a triangular array of numbers 1 to $n$. 
This process can be easily reversed, and is thus a bijection. Monotone triangles can be defined as objects in their own right as follows~\cite{BRESSOUDBOOK}. 

\begin{definition}
Monotone triangles of order $n$
are all triangular arrays of integers with bottom row 1~2~3~\ldots~$n$ 
and integer entries $a_{ij}$ such that $a_{i,j} \le a_{i-1,j} \le a_{i,j+1} \mbox{ and } a_{ij}< a_{i,j+1}$.
\end{definition}

Note that the bottom row of a monotone triangle of order $n$ is always $1~2~3~\ldots~n$. 
If we rotate the monotone triangle clockwise by $\frac{\pi}{4}$ we obtain a semistandard Young tableau (SSYT) of staircase shape $\delta_n = n~(n-1)~(n-2)\ldots 3~2~1$ whose northeast to southwest diagonals are weakly increasing. 
Thus we have the following theorem.

\[
\begin{array}{c}
\mbox{$4\times 4$ ASM}\\
\left( \begin{array}{cccc}
0 & 1 & 0 & 0\\
1 & -1 & 0 & 1\\
0 & 0 & 1 & 0\\
0 & 1 & 0 & 0\end{array} \right)\end{array}
\Longleftrightarrow
\begin{array}{c}
\mbox{Monotone triangle of order 4}\\
\begin{array}{ccccccc}
  & & & 2 & & & \\
  & & 1 & & 4 & & \\
  & 1 & & 3 & & 4 & \\
  1 & & 2 & & 3 & & 4\end{array} \end{array}
\Longleftrightarrow
\begin{array}{c}
\mbox{Rotated array}\\
\begin{array}{cccc}
1&1&1&2\\
2&3&4&\\
3&4&&\\
4&&&
\end{array}\end{array}
\]

\begin{theorem}
\label{prop:asmbij} 
$n\times n$ alternating sign matrices are in bijection with SSYT of staircase shape $\delta_n$ with entries $y_{i,j}$ at most $n$ such that $y_{i,j}\le y_{i+1,j-1}$. Denote this set as $SSA_n$.
\end{theorem}

Ordered by componentwise comparison of the entries, $SSA_n$ forms a distributive lattice $J(P)$ where the Hasse diagram of the
poset of join-irreducibles $P$ (for $n=4$) is shown below:
\begin{center}
\includegraphics[scale=0.25]{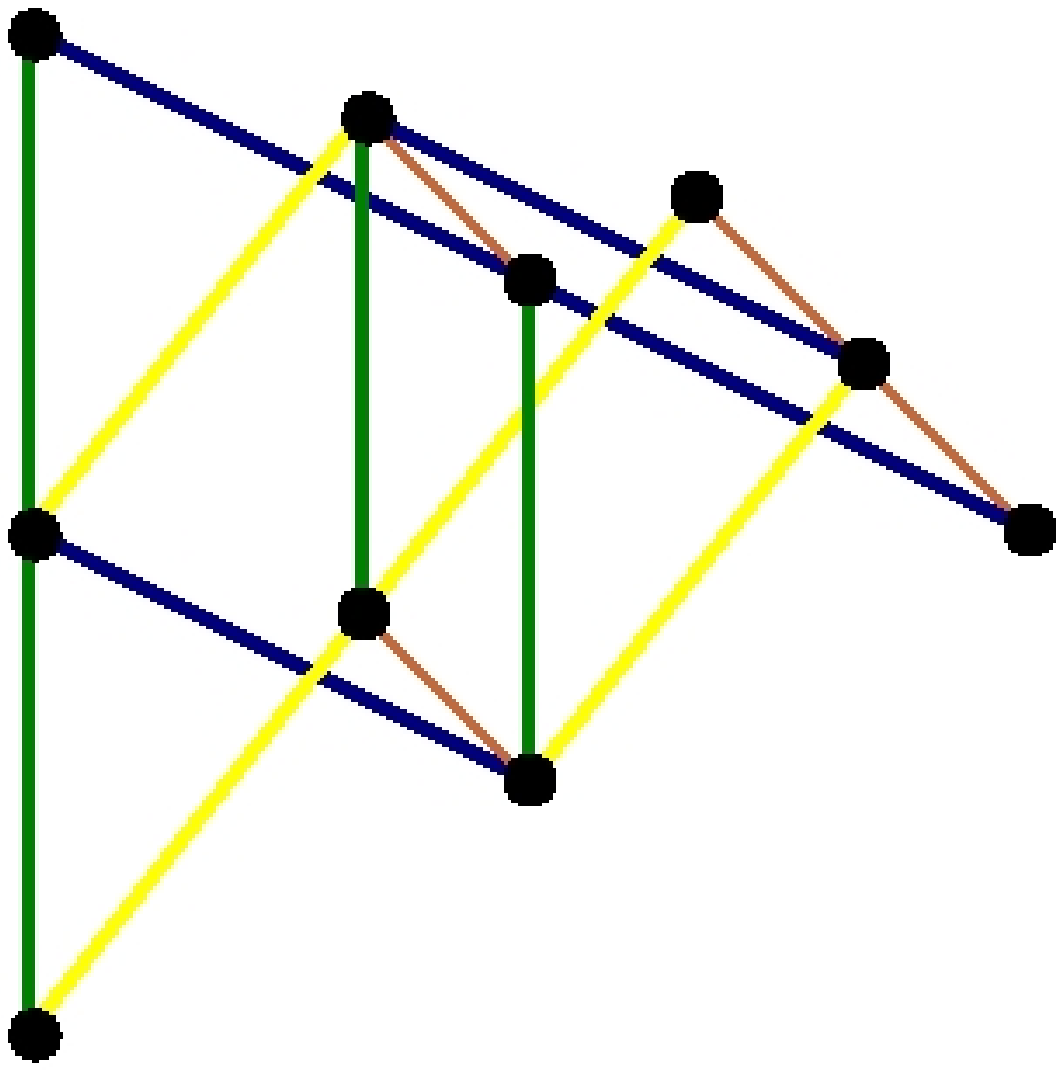}
\end{center}

Given a TSSCPP $t=\{t_{i,j}\}_{1\le i,j\le 2n}$ take a fundamental domain consisting of the triangular array of integers $\{t_{i,j}\}_{n+1\le i\le j\le 2n}$. In this triangular array $t_{i,j}\ge t_{i+1,j}\ge t_{i+1,j+1}$ since $t$ is a plane partition. Also for these values of $i$ and $j$ the entries $t_{i,j}$ satisfy
$0\le t_{i,j}\le 2n+1-i$.  Now if we reflect this array about a vertical line then rotate clockwise by $\frac{\pi}{4}$ we obtain a staircase shape array $x$ whose entries $x_{i,j}$ satisfy the conditions $x_{i,j}\le x_{i,j+1}\le x_{i+1,j}$ and $0\le x_{i,j}\le j$. 
Now add $i$ to each entry in row~$i$. 
This gives us the following theorem.
\begin{theorem}
\label{prop:tsscppbij}
Totally symmetric self--complementary plane partitions inside a $2n\times 2n\times 2n$ box
are in bijection with SSYT of staircase shape $\delta_n$ with entries $y_{i,j}$ at most $n$ such that $y_{i,j}\le y_{i-1,j+1} + 1$. Denote this set as $SST_n$.
\end{theorem}

\[
\begin{array}{c}
\mbox{TSSCPP}\\
\begin{array}{cccccccc}
8&8&8&8&6&6&4&4\\
8&8&8&8&6&5&4&4\\
8&8&7&6&5&4&3&2\\
8&8&6&5&4&3&2&2\\
6&6&5&4&3&2&\cdot&\cdot\\
6&5&4&3&2&1&\cdot&\cdot\\
4&4&3&2&\cdot&\cdot&\cdot&\cdot\\
4&4&2&2&\cdot&\cdot&\cdot&\cdot\\
\end{array}\end{array}
\Longleftrightarrow
\begin{array}{c}
\mbox{Fundamental}\\
\mbox{domain}\\
\begin{array}{cccc}
3&&&\\
2&1&&\\
\cdot&\cdot&\cdot&\\
\cdot&\cdot&\cdot&\cdot\end{array}\end{array}
\Longleftrightarrow
\begin{array}{c}
\mbox{Reflected}\\
\mbox{\& rotated}\\
\begin{array}{cccc}
\cdot&\cdot&1&3\\
\cdot&\cdot&2&\\
\cdot&\cdot&&\\
\cdot&&&
\end{array}
\end{array}
\Longleftrightarrow
\begin{array}{c}
\mbox{$i$ added}\\
\mbox{to row $i$}\\
\begin{array}{cccc}
1&1&2&4\\
2&2&4&\\
3&3&&\\
4&&&
\end{array}\end{array}
\]

Ordered by componentwise comparison of the entries, $SST_n$ forms a distributive lattice $J(Q)$ where the
Hasse diagram of the poset of join-irreducibles $Q$ (for $n=4$) is shown below. (Note that in this paper we will extend the definition of a Hasse diagram slightly by at times drawing edges in the Hasse diagram from $x$ to $y$ when $x<y$ but $y$ does not cover $x$, like the yellow edges below.)

\begin{center}
\includegraphics[scale=0.25]{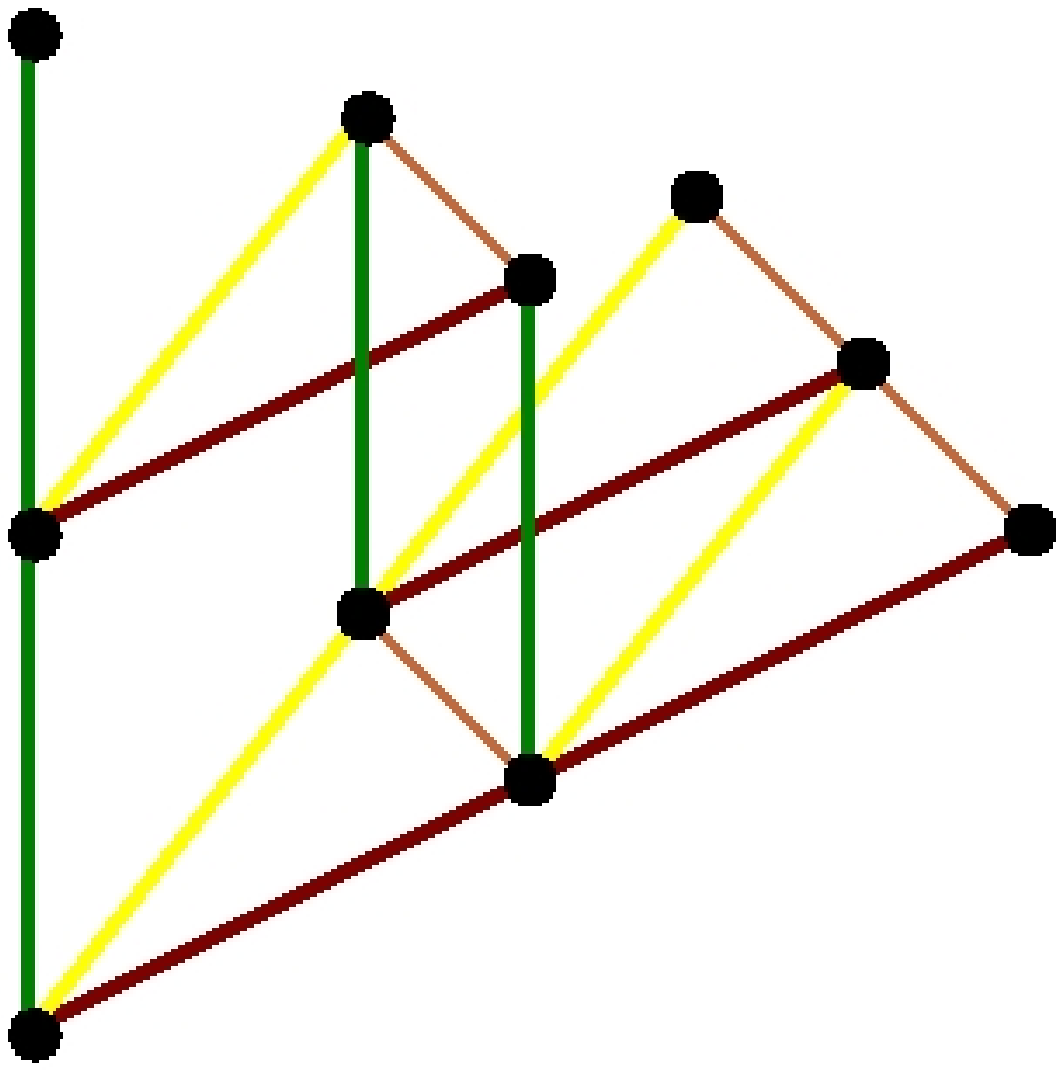}
\end{center}
Suppose we put the posets $P$ and $Q$ together and consider SSYT with both conditions on the diagonals. The Hasse diagram of our new poset looks like a tetrahedron with one direction of edges missing:
\begin{center}
\includegraphics[scale=0.25]{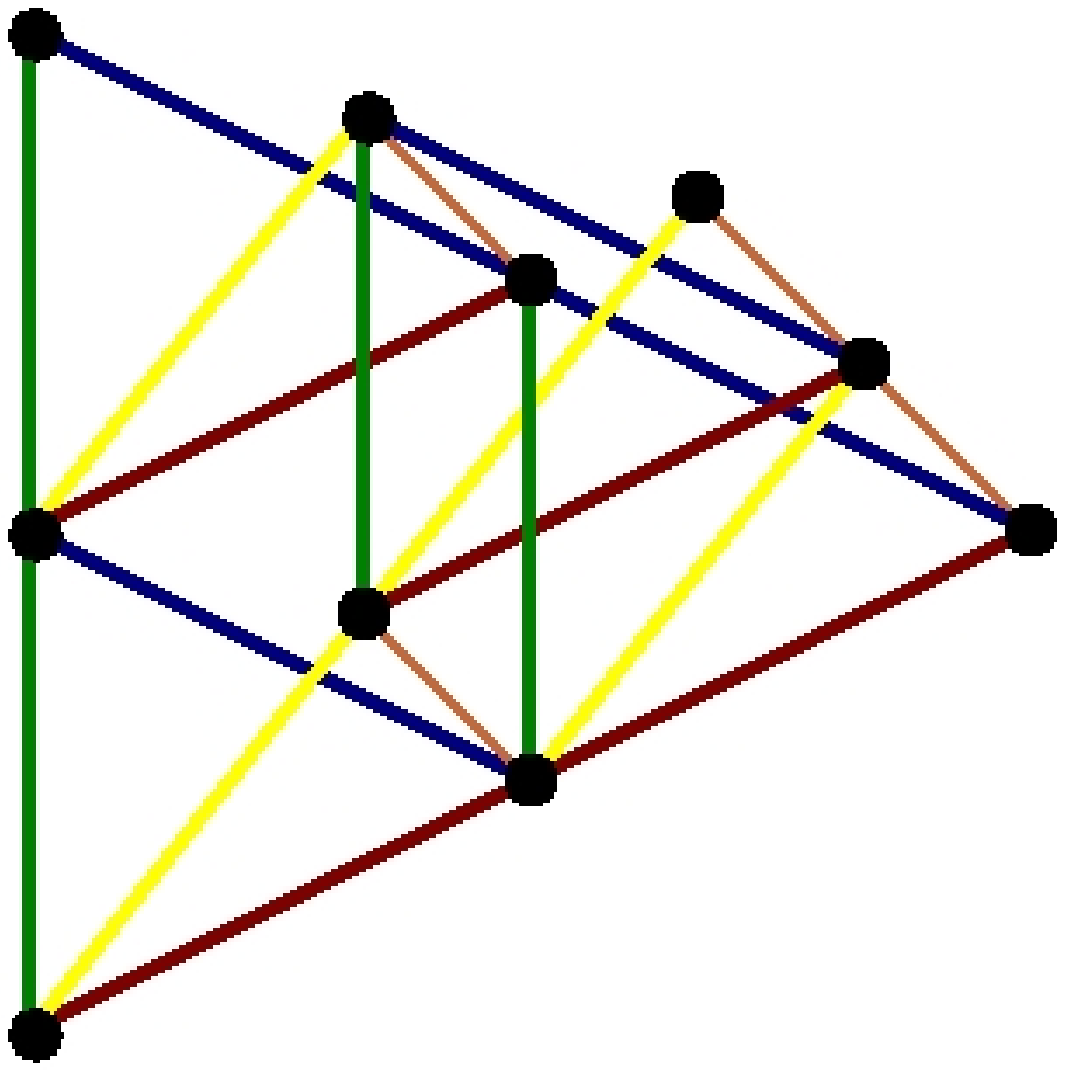}
\end{center}
Inserting those extra edges yields a tetrahedral poset, denoted $T_n$, whose lattice of order ideals we find to be in bijection with totally symmetric plane partitions (TSPPs) inside an $(n - 1) \times (n - 1) \times (n - 1)$ box.

We define $T_n$ precisely as follows. Define the unit vectors $\overrightarrow{r}=(\frac{\sqrt{3}}{2},\frac{1}{2},0)$, $\overrightarrow{g}=(0,1,0)$, and $\overrightarrow{b}=(-\frac{\sqrt{3}}{2},\frac{1}{2},0)$, $\overrightarrow{y}=(\frac{\sqrt{3}}{6},\frac{1}{2},\frac{\sqrt{6}}{3})$, $\overrightarrow{o}=(\frac{-\sqrt{3}}{3},0,\frac{\sqrt{6}}{3})$, $\overrightarrow{s}=(\frac{-\sqrt{3}}{6},\frac{1}{2},-\frac{\sqrt{6}}{3})$. 
Let the elements of $T_n$ be defined as the coordinates of all the points reached by linear combinations of $\overrightarrow{r}$, $\overrightarrow{g}$, and $\overrightarrow{y}$. Thus as a set $T_n =\{c_1 \overrightarrow{r} +c_2 \overrightarrow{g}+c_3\overrightarrow{y}, \mbox{ } c_1,c_2,c_3\in \mathbb{Z}_{\ge 0}, \mbox{ } c_1+c_2+c_3\le n-2\}$. Let all the vectors $\overrightarrow{r}$, $\overrightarrow{g}$, and $\overrightarrow{y}$ used to define the elements of $T_n$ be directed edges in the Hasse diagram of colors red, green, and yellow. Additionally draw as edges of colors blue, orange, and silver the vectors $\overrightarrow{b}$, $\overrightarrow{o}$, and $\overrightarrow{s}$ between poset elements wherever possible. 
The partial order of $T_n$ is defined so that the corner vertex with edges colored red, green, and yellow is the smallest element, the corner vertex with edges colored silver, green, and blue is the largest element, and the other two corner vertices are ordered such that the one with silver, yellow, and orange edges is above the one with orange, red, and blue edges. 
\begin{center}
\includegraphics[scale=0.27]{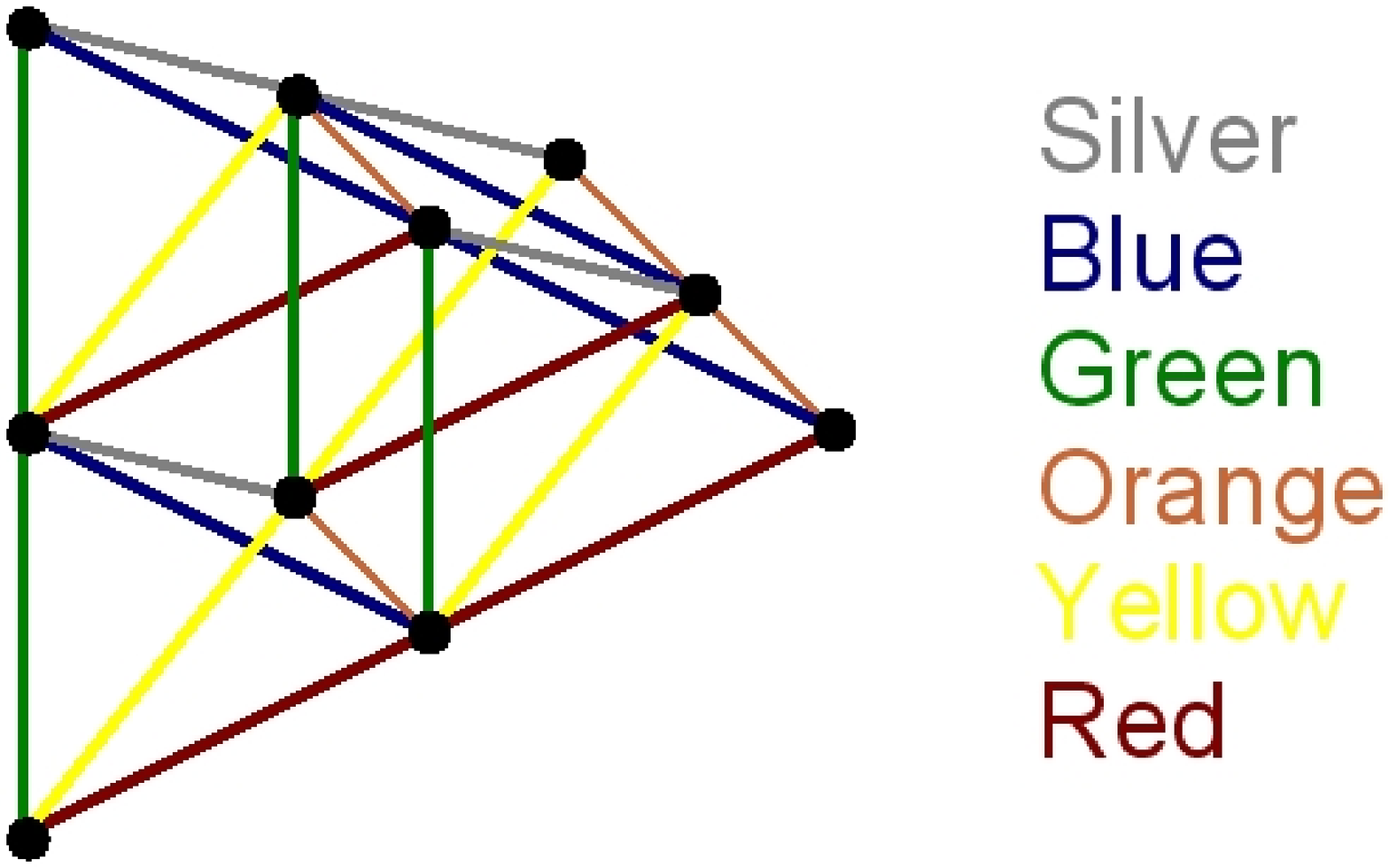}
\end{center}

Since the ASM and TSSCPP posets appear as subposets of $T_n$ with certain edge colors, we now investigate the subposets of $T_n$ made up of all the different combinations of edge colors.
Surprisingly, for almost all of these posets, there exists a nice product formula for the number of order ideals 
and a bijection between these order ideals and an interesting set of combinatorial objects. 
We wish to consider only subsets of the colors which include all the colors whose covering relations are induced by combinations of other colors, 
which are the \emph{admissible} subsets of the following definition.
\begin{definition}
Let a subset $S$ of the six colors $\{$red, blue, green, orange, yellow, silver$\}$ (abbreviated $\{r,b,g,o,y,s\}$) be  \emph{admissible} if all of the following are true: 
If $\{r,b\}\subseteq S$ then $g\in S$, if  $\{o,s\}\subseteq S$ then $b\in S$, if $\{s,y\}\subseteq S$ then $g\in S$, and if $\{r,o\}\subseteq S$ then $y\in S$.
\end{definition}

Given an admissible subset $S$ of the colors $\{r,b,g,o,y,s\}$, let $T_n(S)$ denote the poset formed by the vertices of $T_n$ together with all the edges whose colors are in $S$. The induced colors will be in parentheses.

We give a bijection between order ideals of $T_n(S)$, $S$ an admissible subset of $\{r,b,g,y,o,s\}$, and arrays of integers with certain inequality conditions.
\begin{definition}
\label{def:bofs}
Let $S$ be an admissible subset of $\{r,b,g,y,o,s\}$ and suppose $g\in S$. Define $Y_n(S)$ to be the set of all integer arrays $x$ of staircase shape $\delta_n$ with entries $x_{i,j}$, $1\le i\le n$, $0\le j\le n-i$ satisfying both
$i\le x_{i,j}\le j+i$ and the following inequality conditions corresponding to the additional colors in $S$:
orange: $x_{i,j} < x_{i+1,j}$, red: $x_{i,j}\le x_{i-1,j+1} +1$, yellow: $x_{i,j}\le x_{i,j+1}$, blue: $x_{i,j}\le x_{i+1,j-1}$, silver: $x_{i,j}\le x_{i,j-1}+1$
\end{definition}


For the proof of the following proposition we will need to note the following: $T_n(\{r,b,(g)\})$ is a disjoint poset, whose connected components we will call, from smallest to largest, $P_2, P_3,\ldots P_n$.  Thus $T_n$ can be thought of as the poset which results from beginning with $P_n$, overlaying $P_{n-1},P_{n-2},\ldots,P_3,P_2$ successively, and connecting each $P_{i}$ to $P_{i-1}$ by the orange, yellow, and silver edges.
\begin{proposition}
\label{prop:yplusbij}
If $S$ is an admissible subset of $\{r,b,g,y,o,s\}$ and $g\in S$ then
$Y_n(S)$ 
is in weight--preserving bijection with $J(T_n(S))$ where the weight of $x\in Y_n(S)$ 
is given by $\sum_{i=1}^{n-1} \sum_{j=0}^{n-i} (x_{i,j} - i)$ and the weight of $I\in J(T_n(S))$ equals $|I|$. 
\end{proposition}
\proof
Let $S$ be an admissible subset of $\{r,b,g,y,o,s\}$ and suppose $g\in S$. Recall that $T_n$ is made up of the layers $P_k$ where $2\le k\le n$. Since $g\in S$, $P_k$ contains $k-1$ green--edged chains of length $k-1,\ldots,2,1$. For each $P_k\subseteq T_n$ let the $k-1$ green chains inside $P_k$ determine the entries $x_{i,j}$ ($j\neq 0$) of an integer array on the diagonal where $i+j=k$. In particular, given an order ideal $I$ of $T_n(S)$ form an array $x$ by setting $x_{i,j}$ equal to $i$ plus the number of elements in the induced order ideal of the length~$j$ green chain inside $P_{i+j}$ (in column 0 $x_{i,0}=i$). This defines $x$ as an integer array of staircase shape $\delta_n$ whose entries satisfy $i\le x_{i,j}\le j+i$. Also since each entry $x_{i,j}$ is given by an induced order ideal and since each element of $T_n$ is in exactly one green chain we know that $|I|=\sum_{i,j} x_{i,j}-i$. Thus the weight is preserved.

Now it is left to determine what the other colors mean in terms of the array entries. Since the colors red and blue connect green chains from the same $P_k$ we see that inequalities corresponding to red and blue should relate entries of $x$ on the same northeast to southwest diagonal of $x$. So if $r\in S$ then $x_{i,j}\le x_{i-1,j+1} +1$ and if $b\in S$ then $x_{i,j}\le x_{i+1,j-1}$. The colors yellow, orange, and silver connect $P_k$ to $P_{k+1}$ for $2\le k\le n-1$. So from our construction we see that 
if $o\in S$ then $x_{i,j}\le x_{i+1,j}$, if $y\in S$ then $x_{i,j}\le x_{i,j+1}$, and if $s\in S$ then $x_{i,j}\le x_{i,j-1}+1$.
\endproof

\section{Combinatorial objects as subposet order ideals}

\begin{figure}[htbp]
\centering
\includegraphics[scale=0.45]{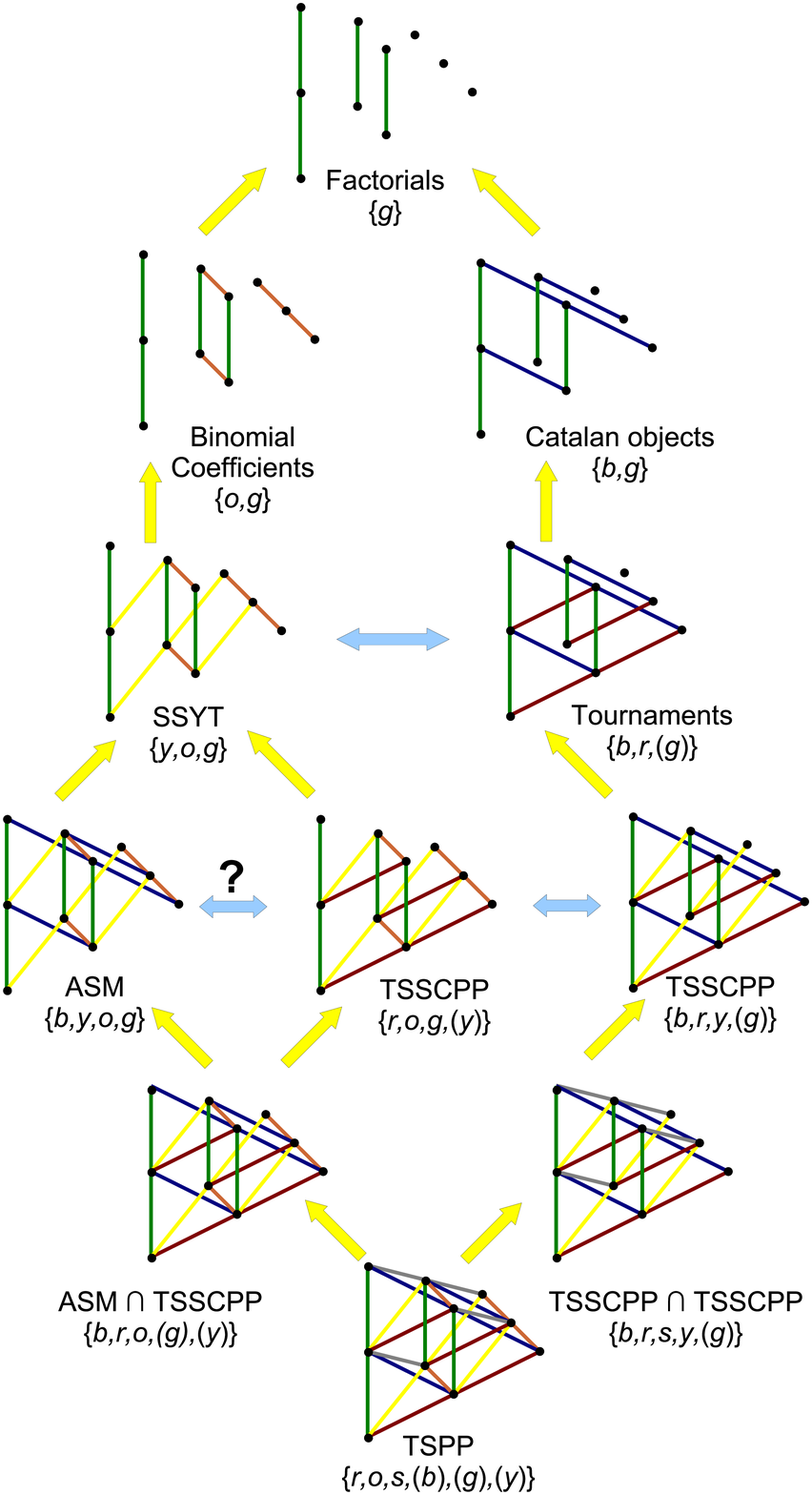}
\caption[The big picture of inclusions and bijections between order ideals]{The big picture of inclusions and bijections between order ideals $J(T_n(S))$. The one sided arrows represent inclusions of one set of order ideals into another. The two sided arrows represent bijections between sets of order ideals. The bijection between the order ideals of the three color posets is in Section~\ref{sec:three} and the bijections between TSSCPP posets is by poset isomorphism. The only missing bijection between sets of order ideals of the same size is between ASM and TSSCPP.}
\label{fig:bigbigpic}
\end{figure}

We will now give product formulas for the number of order ideals of $T_n(S)$ for $S$ an admissible set of colors along with the rank generating functions wherever we have them, where $F(P,q)$ denotes the rank generating function for the poset $P$. For the sake of comparison we have also written each formula as a product over the same indices $1\le i\le j\le k\le n-1$ in a way which is reminiscent of the MacMahon box formula.  See Figure~\ref{fig:bigbigpic} for the big picture of inclusions and bijections between these order ideals. For a more detailed discussion, see~\cite{STRIKER_THESIS}.

\begin{theorem}
For any color $x\in \{r,b,y,g,o,s\}$
\begin{equation}
F(J(T_n(\{x\})),q)=\prod_{j=1}^n j!_q=\displaystyle\prod_{1\le i\le j\le k\le n-1} \frac{[i+1]_q}{[i]_q}.
\end{equation}
\end{theorem}
\proof
$T_n(\{x\})$ is the disjoint sum of $n-j$ chains of length $j-1$ as $j$ goes from 1 to $n-1$. So the number of order ideals is the product of the number of order ideals of each chain. 
\endproof

\begin{theorem}
\label{thm:twoopp}
If $S\in \{\{g,o\},\{r,s\},\{b,y\}\}$ then
\begin{equation}
F(J(T_n(S)),q)=\displaystyle\prod_{j=1}^n {n\brack j}_q=\displaystyle\prod_{1\le i\le j\le k\le n-1} \frac{[j+1]_q}{[j]_q}.
\end{equation}
\end{theorem}
\proof
The arrays $Y_n(\{g,o\})$ strictly decrease down columns and have no conditions on the rows. Thus in a column of length $j$ there must be $j$ distinct integers between $1$ and $n$; this is counted by ${n\choose j}$. If we give a weight to each of these integers of $q$ to the power of that integer minus its row, we have a set $q$--enumerated by the $q$--binomial coefficient ${n\brack j}_q$. Thus $\prod_{j=1}^n {n\brack j}_q$ is the generating function of the arrays $Y_n(\{g,o\})$ and also of the order ideals $F(J(T_n(\{g,o\})),q)$.
The posets $T_n(\{g,o\})$, $T_n(\{r,s\})$, and $T_n(\{b,y\}\}$ are all isomorphic, thus the result follows by poset isomorphism.
\endproof

\begin{theorem}
\label{thm:twoadj}
If $S_1\in \{\{b,g\},\{b,s\},\{y,o\},\{g,s\}\}$ and $S_2\in \{\{r,y\},\{r,g\},\{y,g\},$ $\{b,o\}\}$ then
\begin{equation}
|J(T_n(S_1))|=|J(T_n(S_2))|=\displaystyle\prod_{j=1}^n C_j =\displaystyle\prod_{j=1}^n \frac{1}{j+1} {2j \choose j}=\displaystyle\prod_{1\le i\le j\le k\le n-1} \frac{i+j+2}{i+j}
\end{equation}
\begin{equation}
F(J(T_n(S_1)),q)=F(J(T_n^*(S_2)),q)=\displaystyle\prod_{j=1}^n C_j(q)
\end{equation}
where $^*$ is poset dual, $C_j$ is the $j$th Catalan number, and $C_j(q)$ is the Carlitz--Riordan $q$--Catalan number defined by the recurrence
$C_j(q)=\displaystyle\sum_{k=1}^j q^{k-1} C_{k-1}(q) C_{j-k}(q)$
with initial conditions $C_0(q)=C_1(q)=1$.
\end{theorem}

\proof
$T_n(\{b,g\})$ is isomorphic to the disjoint sum of posets $P_j(\{b,g\})$ for $2\le j\le n$ with rank generating functions $C_j(q)$. 
This can be shown using an easy bijection between these order ideals and Dyck paths from $(0,0)$ to $(2j,0)$.
Thus the number of order ideals $|J(T_n(\{b,g\}))|$ equals the product $\prod_{j=1}^n C_j$ and the rank generating function $F(J(T_n(\{b,g\})),q)$ equals the product $\prod_{j=1}^n C_j(q)$. 
Finally, the posets $T_n(S_1)$ for any choice of $S_1\in \{\{b,g\}$, $\{b,s\}$, $\{y,o\}$, $\{g,s\}\}$ and the posets $T_n^*(S_2)$ for any $S_2\in \{\{r,y\}$, $\{r,g\}$, $\{y,g\}$, $\{b,o\}\}$ are all isomorphic, thus the result follows by poset isomorphism.
\endproof


\begin{theorem}
\label{thm:3nadj}
If $S$ is an admissible subset of $\{r,b,g,o,y,s\}$, $|S|=3$, and $S\notin \{\{r,g,y\},\{s,b,r\}\}$ then
\begin{equation}
F(J(T_n(S)),q)
=\displaystyle\prod_{j=1}^{n-1} (1+q^j)^{n-j}=\displaystyle\prod_{1\le i\le j\le k\le n-1} \frac{[i+j]_q}{[i+j-1]_q}.
\end{equation}
Thus if we set $q=1$ we have $|J(T_n(S))|= 2^{{n\choose 2}}$.
\end{theorem}

We will prove Theorem~\ref{thm:3nadj} using two lemmas since there are two nonisomorphic classes of posets $T_n(S)$ where $S$ is admissible, $|S|=3$, and $S\notin \{\{r,g,y\},\{s,b,r\}\}$. The first lemma is the case where $T_n(S)$ is a disjoint sum of posets and the second lemma is the case where $T_n(S)$ is a connected poset.
\begin{lemma}
\label{lemma:disjoint}
Suppose $S\in \{\{o,s,(b)\}, \{s,y,(g)\}, \{o,r,(y)\}, \{b,r,(g)\}\}$. Then
\[
F(J(T_n(S)),q)
=\displaystyle\prod_{j=1}^{n-1} (1+q^j)^{n-j}=\displaystyle\prod_{1\le i\le j\le k\le n-1} \frac{[i+j]_q}{[i+j-1]_q}.
\]
\end{lemma}
\proof
$T_n(\{b,r,(g)\})$ is a disjoint sum of the $P_j$ posets for $2\le j\le n$. The order ideals of $P_j$ are counted by $2^{j-1}$ and the rank generating function of $J(P_j)$ is given by $\prod_{i=1}^{j-1} (1+q^i)$, both of which are proved by induction. Thus $F(T_n(\{b,r,(g)\}),q)$ is the product of $\prod_{i=1}^{j-1} (1+q^i)$ for $2\le j\le n$. Rewriting the product we obtain $\prod_{j=2}^n\prod_{i=1}^{j-1} (1+q^i)=\prod_{j=1}^{n-1} (1+q^j)^{n-j}$.
The posets $T_n(S)$ where $S\in \{\{o,s,(b)\}, \{s,y,(g)\}, \{o,r,(y)\}\}$ are isomorphic to $T_n(\{b,r,(g)\})$ so the result follows by poset isomorphism.
\endproof

\begin{lemma}
\label{lemma:connected}
Suppose $S\in \{\{r,g,s\},\{o,b,y\},\{y,g,o\},\{b,g,o\},\{y,g,b\}\}$. Then
\[
F(J(T_n(S)),q)
=\displaystyle\prod_{j=1}^{n-1} (1+q^j)^{n-j}=\displaystyle\prod_{1\le i\le j\le k\le n-1} \frac{[i+j]_q}{[i+j-1]_q}.
\]
\end{lemma}
\proof
The arrays $Y_n(\{g,y,o\})$ are by Definition~\ref{def:bofs} equivalent to SSYT of staircase shape $\delta_n$, thus their generating function is given by the Schur function $s_{\delta_n}(x_1,x_2,\ldots,x_n)$. Now
\[
s_{\delta_n}(x_1,x_2,\ldots,x_n)=\frac{\det(x_i^{2(n-j)})_{i,j=1}^n}{\det(x_i^{n-j})_{i,j=1}^n}=\displaystyle\prod_{1\le i<j\le n}\frac{x_i^2-x_j^2}{x_i-x_j}=\displaystyle\prod_{1\le i<j\le n}(x_i+x_j)
\] 
using the algebraic Schur function definition and the Vandermonde determinant. 
The principle specialization of this generating function yields the $q$--generating function $\prod_{j=1}^{n-1} (1+q^j)^{n-j}$.
The posets $T_n(S)$ where $S\in \{\{r,g,s\},\{o,b,y\},\{b,g,o\},\{y,g,b\}\}$ are isomorphic to $T_n(\{g,y,o\})$ so the result follows by poset isomorphism.
\endproof


\proof[Proof of Theorem~\ref{thm:3nadj}]
By Lemma~\ref{lemma:disjoint}, if $S\in \{\{o,s,(b)\}$, $\{s,y,(g)\}$, $\{o,r,(y)\}$, $\{b,r,(g)\}\}$ then the generating function $F(J(T_n(S)),q)$ is as above. By Lemma~\ref{lemma:connected}, if $S\in \{\{r,g,s\}$, $\{o,b,y\}$, $\{y,g,o\}$, $\{b,g,o\}$, $\{y,g,b\}\}$ then generating function $F(J(T_n(S)),q)$ is as above. These are the only admissible subsets $S$ of $\{r,b,g,o,y,s\}$ with $|S|=3$ and $S\notin \{\{r,g,y\},\{s,b,r\}\}$.
\endproof

There seems to be no nice product formula for the number of order ideals of the dual posets $T_n(\{r,g,y\})$ and $T_n(\{s,b,r\})$. 
The number of order ideals up to $n=6$ are: 1, 2, 9, 96, 2498, 161422.

\begin{theorem}
\label{thm:4four}
If $S$ is an admissible subset of $\{r,b,g,o,y,s\}$ and $|S|=4$ then 
\begin{equation}
|J(T_n(S))|=\displaystyle\prod_{j=0}^{n-1} \frac{(3j+1)!}{(n+j)!}=\displaystyle\prod_{1\le i\le j\le k\le n-1} \frac{i+j+k+1}{i+j+k-1}.
\end{equation}
\end{theorem}
\proof
The posets $T_n(S)$ for $S$ admissible and $|S|=4$ are $T_n(\{g,y,b,o\})$, the three isomorphic posets 
$T_n(\{r,o,(y),g\})$, $T_n(\{r,b,(g),y\})$, and $T_n(\{y,s,(g),r\})$, and the three posets dual to these, $T_n(\{y,s,(g),b\})$, $T_n(\{o,s,(b),g\})$, $T_n(\{r,b,(g),s\})$. In Theorem~\ref{prop:asmbij} we showed that the order ideals of $T_n(\{g,y,b,o\})$ are in bijection with $n\times n$ ASMs and in Theorem~\ref{prop:tsscppbij} we showed that the order ideals of $T_n(\{r,o,(y),g\})$ are in bijection with TSSCPPs inside a $2n\times 2n\times 2n$ box.  Therefore by poset isomorphism TSSCPPs inside a $2n\times 2n\times 2n$ box are in bijection with the order ideals of any of  $T_n(\{r,o,(y),g\})$, $T_n(\{r,b,(g),y\})$, $T_n(\{y,s,(g),r\})$, $T_n^*(\{y,s,(g),b\})$, $T_n^*(\{o,s,(b),g\})$, or  $T_n^*(\{r,b,(g),s\})$. Thus by the enumeration of ASMs in~\cite{ZEILASM} and~\cite{KUP_ASM_CONJ} and the enumeration of TSSCPPs in~\cite{ANDREWS_PPV} we have the above formula for the number of order ideals.
\endproof

There are two different cases for five colors: one case consists of the dual posets $T_n(\{(g),(b),o,y,s\})$ and $T_n(\{r,b,(g),o,(y)\})$ and the other case is $T_n(\{r,b,s,(y),g\})$. A nice product formula has not yet been found for either case. 

\begin{theorem}
\label{thm:6six}
\begin{equation}
|J(T_n)|
=\displaystyle\prod_{1\le i\le j\le n-1} \frac{i+j+n-2}{i+2j-2}=\displaystyle\prod_{1\le i\le j\le k\le n-1} \frac{i+j+k-1}{i+j+k-2}
\end{equation}
\end{theorem}

\proof
Totally symmetric plane partitions are plane partitions which are symmetric with respect to all permutations of the $x,y,z$ axes. Thus we can take as a fundamental domain the wedge where $x\ge y\ge z$. Then if we draw the lattice points in this wedge (inside a fixed bounding box of size $n-1$) as a poset with edges in the $x$, $y$, and $z$ directions, we obtain the poset $T_n$ where the $x$ direction corresponds to the red edges of $T_n$, the $y$ direction to the orange edges, and the $z$ direction to the silver edges. All other colors of edges in $T_n$ are induced by the colors red, silver, and orange. Thus TSPPs inside an $(n-1)\times (n-1)\times (n-1)$ box are in bijection with the order ideals of $T_n$.
Thus by the enumeration of TSPPs in~\cite{STEMBRIDGE_TSPP} the number of order ideals $|J(T_n)|$ is given by the above formula.
\endproof

\section{Bijections with tournaments}
\label{sec:three}
Theorem~\ref{thm:3nadj} states that the order ideals of the three color posets $J(T_n(S))$ where $S$ is admissible, $|S|=3$, and $S\notin \{\{r,g,y\},\{s,b,r\}\}$ are counted by 
$2^{n\choose 2}$. This is also the number of graphs on $n$ labeled vertices and equivalently the number of tournaments on $n$ vertices. A tournament is a complete directed graph with labeled vertices. 
We now discuss bijections between these order ideals 
and tournaments. 
\begin{theorem}
\label{thm:easy2nc2bij}
There exists an explicit 
bijection between the order ideals of the poset $T_n(\{b,r,(g)\})$ and tournaments on $n$ vertices. 
\end{theorem}
\[
\begin{array}{ccc}
1&1 &1 \\
2&2& \\
3&&
\end{array}
\hspace{.3cm}
\begin{array}{ccc}
1&1 &\textcolor{red}{2} \\
2&2&\\
3&&
\end{array}
\hspace{.3cm}
\begin{array}{ccc}
1&1 &2 \\
2&\textcolor{red}{3}&\\
3&&
\end{array}
\hspace{.3cm}
\begin{array}{ccc}
1&1 & \textcolor{red}{3} \\
2&\textcolor{red}{3}&\\
3&&
\end{array}
\hspace{.3cm}
\begin{array}{ccc}
1&\textcolor{red}{2} & 1 \\
2&2&\\
3&&
\end{array}
\hspace{.3cm}
\begin{array}{ccc}
1&\textcolor{red}{2} & \textcolor{red}{2} \\
2&2&\\
3&&
\end{array}
\hspace{.3cm}
\begin{array}{ccc}
1&\textcolor{red}{2} & 2 \\
2&\textcolor{red}{3}&\\
3&&
\end{array}
\hspace{.3cm}
\begin{array}{ccc}
1&\textcolor{red}{2} & \textcolor{red}{3} \\
2&\textcolor{red}{3}&\\
3&&
\end{array}
\]

\begin{center}
\includegraphics[scale=0.28]{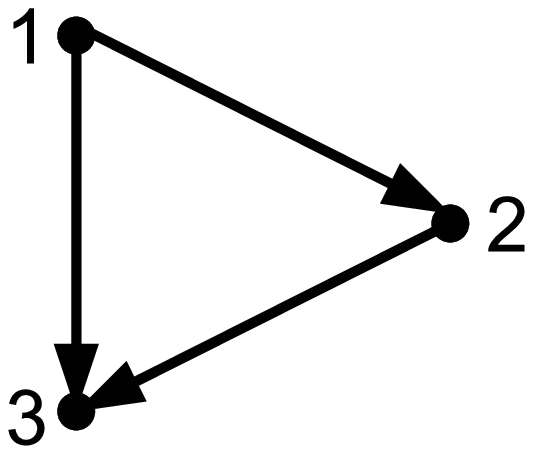}
\includegraphics[scale=0.28]{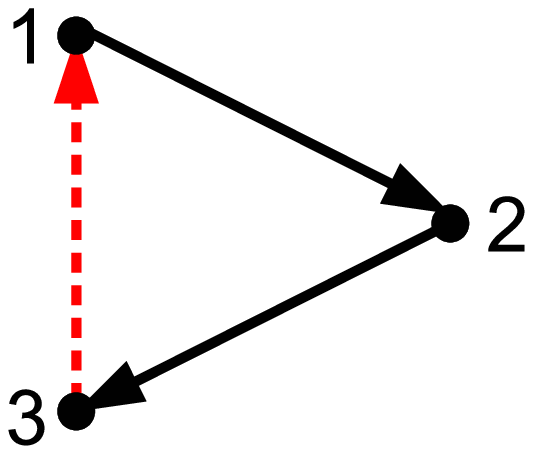}
\includegraphics[scale=0.28]{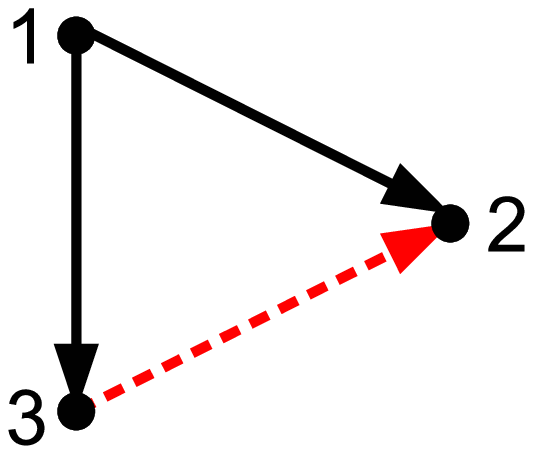}
\includegraphics[scale=0.28]{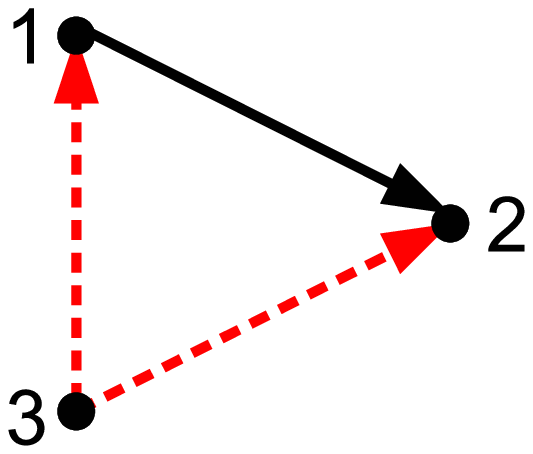}
\includegraphics[scale=0.28]{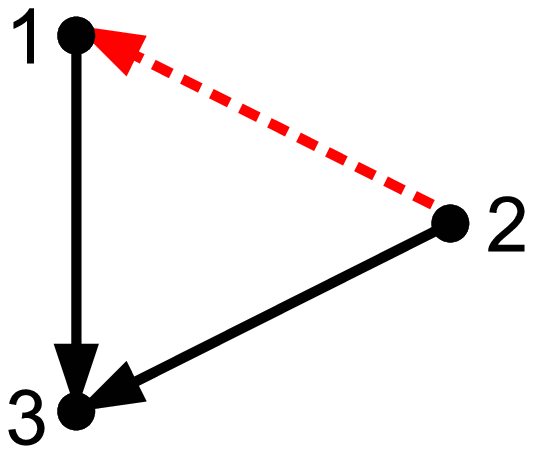}
\includegraphics[scale=0.28]{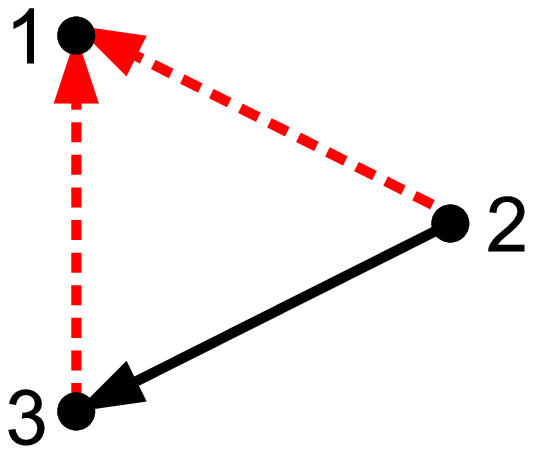}
\includegraphics[scale=0.28]{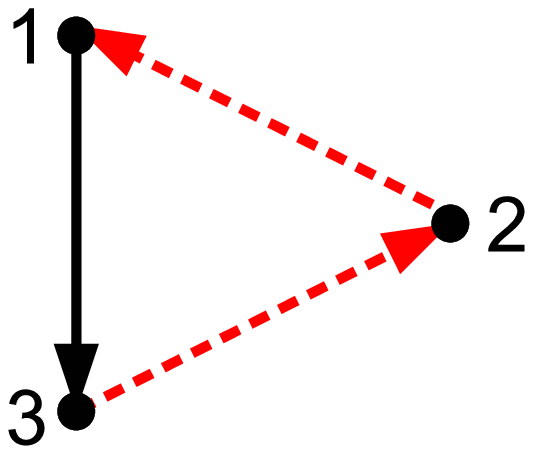}
\includegraphics[scale=0.28]{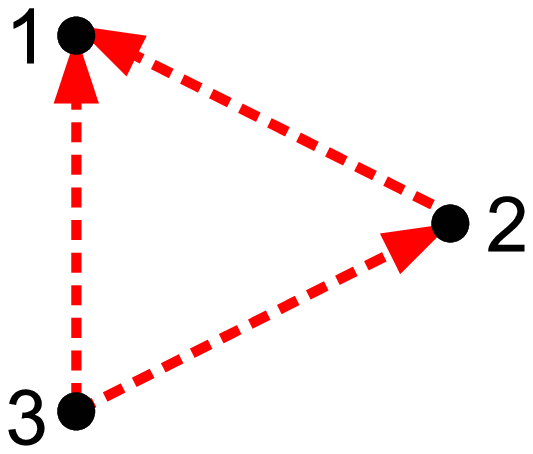}
\end{center}
\proof 
The colors blue and red correspond to inequalities on $Y_n(\{b,r,(g)\})$ such that as one goes up the southwest to northeast diagonals at each step the next entry has the choice between staying the same and decreasing by one. 
Therefore since each of the ${n\choose 2}$ entries of the array not in the 0th column has exactly two choices of values given the value of the entry to the southwest, we may consider each array entry $\alpha_{i,j}$ with $j\ge 1$ to symbolize the outcome of the game between $i$ and $i+j$ in a tournament. If $\alpha_{i,j} = \alpha_{i+1,j-1}$ say the outcome of the game between $i$ and $i+j$ is an upset and otherwise not. Thus tournaments on $n$ vertices are in bijection with the arrays $Y_n(\{b,r,(g)\})$ and also with the order ideals of $T_n(\{b,r,(g)\})$. 
\endproof

The bijection between the order ideals of the connected three color posets of Lemma~\ref{lemma:connected} and tournaments is due to Sundquist in~\cite{SUNDQUIST} and involves repeated use of jeu de taquin and column deletion to go from SSYT of shape $\delta_n$ and largest entry $n$ to certain tableaux in bijection with tournaments on $n$ vertices. 

Next we describe which subsets of tournaments correspond to TSSCPPs.
\begin{theorem}
\label{thm:tsscpptournbij}
TSSCPPs inside a $2n\times 2n\times 2n$ box are in bijection with tournaments on vertices labeled $1,2,\ldots,n$ which satisfy the following condition on the upsets: if vertex $v$ has $k$ upsets with vertices in $\{u,u+1,\ldots,v-1\}$ then vertex $v-1$ has at most $k$ upsets with vertices in $\{u,u+1,\ldots,v-2\}$. 
\end{theorem}
\proof
We have seen in Theorem~\ref{thm:easy2nc2bij} the bijection between the order ideals of $T_n(\{r,b,(g)\})$ and tournaments on $n$ vertices. 
Thus if we consider the TSSCPP arrays $Y_n(\{r,b,(g),y\})$ we need only find an interpretation for the yellow edges in terms of tournaments. 
Recall that yellow corresponds to a weak increase across the rows of $\alpha$. To satisfy this condition, for each choice of $i\in\{1,\ldots,n-1\}$ and $j\in\{1,\ldots,n-i-1\}$ the number of diagonal equalities to the southwest of $\alpha_{i,j}$ must be less than or equal to the number of diagonal equalities to the southwest of $\alpha_{i,j+1}$. So in terms of tournaments, the number of upsets between $i+j+1$ and vertices greater than or equal to $i$ must be greater than or equal to the number of upsets between $i+j$ and vertices greater than or equal to $i$.
\endproof

\section{Connections between ASMs, TSSCPPs, and tournaments}
\label{sec:3to4}

In this section we discuss the expansion of the tournament generating function as a sum over ASMs and derive a new expansion as a sum over TSSCPPs. 
We begin with the following theorem of Robbins and Rumsey~\cite{ROBBINSRUMSEY}. We need the following notion: the inversion number of an ASM $A$ is defined as \linebreak $I(A)=\sum A_{ij} A_{k\ell}$ where the sum is over all $i,j,k,\ell$ such that $i>k$ and $j<\ell$. 

\begin{theorem}[Robbins--Rumsey]
\label{thm:lambdadet}
Let $A_n$ be the set of $n\times n$ alternating sign matrices, and for $A\in A_n$ let $I(A)$ denote the inversion number of $A$ and $N(A)$ the number of $-1$ entries in $A$, then
\begin{equation}
\label{eq:2enumeq}
\displaystyle\prod_{1\le i < j\le n} (x_i + \lambda x_j) = \displaystyle\sum_{A\in A_n} \lambda^{I(A)}\left(1+\lambda^{-1}\right)^{N(A)}\displaystyle\prod_{i,j=1}^n x_{j}^{(n-i)A_{ij}}.
\end{equation}
\end{theorem}

Note that 
the left--hand side is the generating function for tournaments on $n$ vertices where each factor of $(x_i + \lambda x_j)$ represents the outcome of the game between $i$ and $j$ in the tournament. If $x_i$ is chosen then the expected winner, $i$, is the actual winner, and if $\lambda x_j$ is chosen then $j$ is the unexpected winner and the game is an upset. Thus in each monomial in the expansion of $\prod_{1\le i < j\le n}(x_i + \lambda x_j)$ the power of $\lambda$ equals the number of upsets and the power of $x_k$ equals the number of wins of $k$. 

We rewrite Theorem~\ref{thm:lambdadet} in different notation which will also be needed later.
For any staircase shape integer array $\alpha\in Y_n(S)$ 
let $E_{i,k}(\alpha)$ be the number of entries of value $k$ in row $i$ equal to their southwest diagonal neighbor, $E^i(\alpha)$ be the number of entries in (southwest to northeast) diagonal $i$ equal to their southwest diagonal neighbor, and $E_i(\alpha)$ be the number of entries in row $i$ equal to their southwest diagonal neighbor, that is, $E_i(\alpha)=\sum_k E_{i,k}(\alpha)$. Also let $E(\alpha)$ be the total number of entries of $\alpha$ equal to their southwest diagonal neighbor, that is, $E(\alpha)=\sum_{i} E_{i}(\alpha)=\sum_{i} E^{i}(\alpha)$. We now define variables for the content of $\alpha$. Let $C_{i,k}(\alpha)$ be the number of entries in row $i$ with value $k$ and let $C_k(\alpha)$ be the total number of entries of $\alpha$ equal to $k$, that is, $C_{k}(\alpha)=\sum_i C_{i,k}(\alpha)$. Let $N(\alpha)$ be the number of entries of $\alpha$ strictly greater than their neighbor to the west and strictly less than their neighbor to the southwest. When $\alpha\in Y_n(\{b,y,o,g\})$ then $N(\alpha)$ equals the number of $-1$ entries in the corresponding ASM. 
%
%

Using this notation we reformulate Theorem~\ref{thm:lambdadet} in the following way.
\begin{theorem}
\label{thm:lambdaASM}
The generating function for tournaments on $n$ vertices can be expanded as a sum over the ASM arrays $Y_n(\{b,y,o,g\})$ in the following way.
\begin{equation}
\label{eq:asmlambda}
\displaystyle\prod_{1\le i < j\le n} (x_i + \lambda x_j) = \displaystyle\sum_{\alpha\in Y_n(\{b,y,o,g\})} \lambda^{E(\alpha)} (1 + \lambda)^{N(\alpha)} \displaystyle\prod_{k=1}^n x_k^{C_k(\alpha)-1}
\end{equation}
\end{theorem}
\proof
First we rewrite Equation~\ref{eq:2enumeq} by factoring out $\lambda^{-1}$ from each $\left(1+\lambda^{-1}\right)$.
\[
\displaystyle\prod_{1\le i < j\le n} (x_i + \lambda x_j) = \displaystyle\sum_{A\in A_n} \lambda^{I(A) - N(A)}\left(1+\lambda\right)^{N(A)}\displaystyle\prod_{i,j=1}^n x_{j}^{(n-i)A_{ij}}
\]
Let $\alpha\in Y_n(\{b,y,o,g\})$ be the array which corresponds to $A$. It is left to show that $I(A)-N(A)=E(\alpha)$ and $\prod_{i,j=1}^n x_{j}^{(n-i)A_{ij}}=\prod_{j=1}^n x_j^{C_j(\alpha)-1}$. In the latter equality take the product over $i$ of the left hand side: $\prod_{i,j=1}^n x_{j}^{(n-i)A_{ij}}=\prod_{j=1}^n x_{j}^{\sum_{i=1}^n(n-i)A_{ij}}$. 
We wish to show $C_{j}(\alpha)-1=\sum_{i=1}^n (n-i)A_{ij}$. $C_j(\alpha)$ equals the number of entries of $\alpha$ with value $j$, so $C_j(\alpha)-1$ equals the number of entries of $\alpha$ with value $j$ not counting the $j$ in the 0th column. 
Now
the number of $j$s in columns 1 through $n-1$ of $\alpha$ equals the number of 1s in column $j$ of $A$ plus the number of zeros in column $j$ of $A$ which are south of a 1 with no $-1$s in between. This is precisely what $\sum_{i=1}^n(n-i)A_{ij}$ counts by taking a positive contribution from every~1 and every entry below that~1 in column $j$ and then subtracting one for every~$-1$ and every entry below that $-1$ in column $j$. Thus $C_{j}(\alpha)-1=\sum_{i=1}^n (n-i)A_{ij}$ so that $\prod_{i,j=1}^n x_{j}^{(n-i)A_{ij}}=\prod_{j=1}^n x_j^{C_j(\alpha)-1}$.

Now we wish to show that $I(A)-N(A)=E(\alpha)$. 
Fix $i$, $j$, and $\ell$ and consider $\sum_{k<i} A_{ij} A_{k\ell}$. Let $k'$ be the row of the southernmost nonzero entry in column $\ell$ such that $k'<i$. If there exists no such $k'$ (that is, $A_{k\ell}=0 \mbox{ }\forall\mbox{ } k<i$) or if $A_{k'\ell}=-1$ then $\sum_{k>i} A_{ij} A_{k\ell}=0$ since there must be an even number of nonzero entries in $\{A_{k\ell}, k<i\}$ half of which are 1 and half of which are $-1$. If $A_{k'\ell}=1$ then $\sum_{k<i} A_{ij} A_{k\ell}=A_{ij}$. Thus $I(A)=\sum_{i,j} \alpha_{ij} A_{ij}$ where $\alpha_{ij}$ equals the number of columns east of column $j$ such that $A_{k'\ell}$ with $k'>i$ exists and equals~1. Let column $\ell'$ be one of the columns counted by $\alpha_{ij}$. Then $A_{i\ell'}$ cannot equal~1, otherwise $A_{k'\ell'}$ would either not exist or equal $-1$. If $A_{i\ell'}=0$ then in $\alpha$ there is a corresponding diagonal equality. If $A_{i\ell'}=-1$ then there is no diagonal equality in $\alpha$. Thus $I(A)=E(\alpha)+N(A)$.
\endproof


Many people have wondered what the TSSCPP analogue of the $-1$ in an ASM may be. The following theorem does not give a direct analogue, but rather expands the left--hand side of~(\ref{eq:asmlambda}) as a sum over TSSCPPs instead of ASMs.

\begin{theorem}
\label{thm:lambdaTSSCPP}
The generating function for tournaments on $n$ vertices can be expanded as a sum over the TSSCPP arrays $Y_n(\{b,r,(g),y\})$ in the following way.
\begin{equation}
\label{eq:lambdatsscpp}
\displaystyle\prod_{1\le i < j\le n} (x_i + \lambda x_j) = \displaystyle\sum_{\alpha\in Y_n(\{b,r,(g),y\})} \lambda^{E(\alpha)} \displaystyle\prod_{i=1}^{n-1} x_i^{n-i-E_i(\alpha)} \displaystyle\sum_{\mbox{row shuffles $\alpha'$ of $\alpha$}} \mbox{ }\displaystyle\prod_{j=1}^{n-1} x_j^{E^j(\alpha')}
\end{equation}
where a row shuffle $\alpha'$ of $\alpha\in Y_n(\{b,r,(g),y\}$ is an array obtained by reordering the entries in the rows of $\alpha$ in such a way that $\alpha'\in Y_n(\{b,r,(g)\}$. Also, setting the $x$'s to 1 we have 
\begin{equation}
\label{eq:lambdanoxtsscpp}
(1+\lambda)^{n\choose 2} = \displaystyle\sum_{\alpha\in Y_n(\{b,r,(g),y\})} \lambda^{E(\alpha)} \displaystyle\prod_{1\le i\le k\le n-1} {C_{i+1,k}(\alpha)\choose E_{i,k}(\alpha)}. 
\end{equation}
\end{theorem}
\proof
We begin with the set $Y_n(\{b,r,(g),y\})$ and remove the inequality restriction corresponding to the color yellow to obtain the arrays $Y_n(\{b,r,(g)\})$ (which are in bijection with tournaments). 
We use the following algorithm for turning any $\alpha\in Y_n(\{b,r,(g)\})$ into an element of $Y_n(\{b,r,(g),y\})$ thus grouping all the elements of $Y_n(\{b,r,(g)\})$ into fibers over the elements of $Y_n(\{b,r,(g),y\})$. Assume each row of $\alpha$ below row $i$ is weakly increasing. Thus $\alpha_{i+1,j}\le \alpha_{i+1,j+1}$. If $\alpha_{i+1,j}< \alpha_{i+1,j+1}$ then $\alpha_{i,j+1}\le \alpha_{i+1,j+2}$ since $\alpha_{i,j+1}\in \{\alpha_{i+1,j}, \alpha_{i+1,j}-1\}$ and $\alpha_{i,j+2}\in \{\alpha_{i+1,j+1},\alpha_{i+1,j+1}-1\}$ by the inequalities corresponding to red and blue. So the only entries which may be out of order in row~$i$ are those for which their southwest neighbors are equal. If $\alpha_{i+1,j} = \alpha_{i+1,j+1}$ but $\alpha_{i,j+1}>\alpha_{i,j+2}$ it must be that $\alpha_{i,j+1}=\alpha_{i+1,j}$ and $\alpha_{i,j+2}=\alpha_{i+1,j+1}-1$. So we may swap $\alpha_{i,j+1}$ and $\alpha_{i,j+2}$ along with their entire northeast diagonals while not violating the red and blue inequalities. 
By completing this process for all rows we obtain an array with weakly increasing rows which is thus in $Y_n(\{b,r,(g),y\})$. 

Now we do a weighted count of how many arrays in $Y_n(\{b,r,(g)\})$ are mapped to a given array in $Y_n(\{b,r,(g),y\})$. Again we rely on the fact that entries in a row can be reordered only when their southwest neighbors are equal. Thus to find the weight of all the $Y_n(\{b,r,(g)\})$ arrays corresponding to a single $Y_n(\{b,r,(g),y\})$ array we simply need to find the set of diagonals containing equalities. The diagonal equalities give a weight dependent on which diagonal they are in, whereas the diagonal inequalities give a weight according to their row (which remains constant). Thus if we are keeping track of the $x_i$ weight we can do no better than to write this as a sum over all the allowable 
shuffles of the rows of $\alpha$ with the $x$ weight of the diagonal equalities dependent on the position. Thus we have Equation (\ref{eq:lambdatsscpp}). 

If we set $x_i=1$ for all $i$ and only keep track of the $\lambda$ we can make a more precise statement. The above proof shows that the $\lambda$'s result from the diagonal equalities, and the number of different reorderings of the rows tell us the number of different elements of $Y_n(\{b,r,(g)\})$ which correspond to a given element of $Y_n(\{b,r,(g),y\})$. We count this number of allowable reorderings as a product over all rows $i$ and all array values $k$ as ${C_{i+1,k}(\alpha)\choose E_{i,k}(\alpha)}$. This yields Equation~(\ref{eq:lambdanoxtsscpp}).
\endproof


The difference in the weighting of ASMs and TSSCPPs in Theorems~\ref{thm:lambdaASM} and~\ref{thm:lambdaTSSCPP} is substantial. For ASMs the more complicated part of the formula arises in the power of $\lambda$ and for TSSCPPs the complication comes from the $x$ variables. These theorems are also strangely similar. They show that the tournament generating function can be expanded as a sum over either ASMs or TSSCPPs, but we still have no direct reason why the number of summands should be the same. 
The combination of Theorems~\ref{thm:lambdaASM} and \ref{thm:lambdaTSSCPP} may contribute toward finding a bijection between ASMs and TSSCPPs, but the  differences between these expansions show why a bijection is not obvious.




\bibliographystyle{amsalpha}

\end{document}